\documentclass[12pt]{article}

\usepackage{amsmath,amssymb}
\usepackage[applemac]{inputenc}
\usepackage{amsmath,amssymb,fullpage}
\usepackage{showlabels}
\usepackage{color}

\newcommand{\dproof}{\noindent {Proof.} \quad}
\newcommand{\fproof}{\hfill $\square$ \bigskip}

\newtheorem{definition}{Definition}[section]

\newtheorem{theorem}[definition]{Theorem}
\newtheorem{problem}[definition]{Problem}

\numberwithin{equation}{section}

\def\1B{\text{1\!\!I}}

\begin{document}
\date{26 March 2017 }
\title{Optimal insider control of stochastic Volterra equations}

\author{
Olfa Draouil$^{1}$ }

\footnotetext[1]{Department of Mathematics, University of Tunis El Manar, Tunis, Tunisia.\\
Email: {\tt olfadraouil@hotmail.fr}}



\maketitle

\paragraph{MSC(2010):} 60H10, 91A15, 91A23, 91B38, 91B55, 91B70, 93E20

\paragraph{Keywords:} Stochastic Volterra equation; optimal control; inside information; Donsker delta functional; stochastic maximum principle.

\begin{abstract} We study the problem of optimal inside control of a stochastic Volterra equation driven by a Brownian motion and a Poisson random measure.
We prove a sufficient and a necessary maximum principle for the optimal control when the trader has only partial information available to her decisions and on the other hand, may have some inside information about the
future of the system.
The results are applied to the problem of finding the optimal insider portfolio in a financial market where the risky asset price is given by a stochastic Volterra equation.

\end{abstract}

\section{Introduction}

In this paper we consider an optimal control problem for a stochastic process $X(t)=X(t,Z)=X(t,z)|_{z=Z}$ defined as the solution of a stochastic Volterra (integral) equation (SVIE) given by
\begin{align}\label{eq1.1a}
&X(t) = X^u(t,Z) = \xi(t) +\int_0^ t b (t, s, X(s,Z), u(s,Z)) ds\nonumber\\
&+ \int_0^ t \sigma (t, s, X(s,Z), u(s,Z)) dB(s) + \int_0^ t\int_{\mathbb{R}} \gamma(t, s, X(s,Z), u(s,Z), \zeta) \tilde{N}(ds, d\zeta)
\end{align}
 Here $B(t)$ and $\tilde{N}(dt, d\zeta)$ is a Brownian motion and an independent compensated Poisson random measure, respectively, jointly defined on a filtered probability space $(\Omega, \mathbb{F} =\{\mathcal{F}_t\}_{t\geq0}, P)$ satisfying the usual conditions.
 The process $u(t) = u(t,Z) = u(t, z)_{z=Z}$ is our insider control
process, where $Z$ is a given $\mathcal{F}_{T_0}$-measurable random variable for some $T_0 > 0$ , representing
the inside information available to the controller.
Note that from \eqref{eq1.1a} we get:
\begin{align}
&dX(t) = \xi'(t)dt + b (t, t, X(t,Z), u(t,Z)) dt + (\int_0^t \frac{\partial b}{\partial t} (t, s, X(s,Z), u(s,Z)) ds) dt\nonumber\\
&+\sigma(t, t, X(t,Z), u(t,Z)) dB(t) + (\int_0^t \frac{\partial \sigma}{\partial t}(t, s, X(s,Z), u(s,Z)) dB(s)) dt\nonumber\\
&+\int_{\mathbb{R}} \gamma(t, t, X(t,Z), u(t,Z),\zeta) \tilde{N}(dt, d\zeta) + (\int_0^t\int_{\mathbb{R}} \frac{\partial \gamma}{\partial t} (t, s, X(s,Z), u(s,Z),\zeta) \tilde{N}(ds, d\zeta)) dt.
\end{align}

We assume that the inside information is of \emph{initial enlargement} type. Specifically, we assume  that the inside filtration $\mathbb{H}$ has the form

\begin{equation}\label{eq1.1}
 \mathbb{H}= \{ \mathcal{H}_t\}_{0\leq t \leq T}, \text{ where } \mathcal{H}_t = \mathcal{F}_t \vee \sigma(Z)
\end{equation}
for all $t$, where $Z$ is a given $\mathcal{F}_{T_0}$-measurable random variable, for some $T_0 > 0$ (constant).
Here and in the following we use the right-continuous version of  $ \mathbb{H}$, i.e. we put
$\mathcal{H}_{t}= \mathcal{H}_{t^+}=\bigcap_{s>t}\mathcal{H}_s.$

We also assume that the \emph{Donsker delta functional} of $Z$ exists (see below). This assumption implies that the Jacod condition holds, and hence that $B(\cdot)$ and $N(\cdot,\cdot)$ are semimartingales with respect to $\mathbb{H}$. See e.g. \cite{DO3} for details.
Let $\mathcal{A}$ be a given family of admissible controls, required to be $\mathcal{K}_t=\mathcal{G}_t\vee\sigma(Z)$-predictable, where $\mathbb{G}=\{\mathcal{G}_t\}_{0\leq t\leq T}$ is a given subfiltration of $\mathbb{F}=\{\mathcal{F}_t\}_{0\leq t\leq T}$ in the sense that $\mathcal{G}_t\subseteq \mathcal{F}_t$ for all $t \in[0,T]$. That is mean that $\mathcal{K}_t\subseteq \mathcal{H}_t$ for all $t\in[0,T]$.
We assume that the value at time $t$ of our insider control process $u(t,x)$ is allowed to depend on both $Z$ and $\mathcal{G}_t$. In other words, $u(.,x)$ is assumed to be $\mathbb{K}$-adapted. Therefore it has the form
\begin{equation}\label{eq1.2}
    u(t,\omega) = u_1(t, Z, \omega)
\end{equation}
for some function $u_1 : [0, T]\times \mathbb{R}\times\Omega\rightarrow \mathbb{R}$
such that $u_1(. , z)$ is $\mathbb{G}$-adapted for each $z \in\mathbb{R}$. For simplicity (albeit with some abuse of notation) we will in the following write $u$ instead of $u_1$.
In other word, we study the case when the controller  has only partial information available to her decisions and on the
other hand, may have some inside information about the
future of the system.

Let $\mathbb{U}$ denote the set of admissible control values. We assume that the functions
\begin{align}
b(t,s,x,u,z)&=b(t,s,x,u,z,\omega):
[0,T]\times [0,T] \times \mathbb{R} \times \mathbb{U}\times \mathbb{R} \times \Omega \mapsto \mathbb{R}\nonumber\\
\sigma(t,s,x,u,z)&=\sigma(t,s,x,u,z,\omega):
[0,T]\times [0,T]\times \mathbb{R} \times \mathbb{U} \times \mathbb{R} \times \Omega \mapsto \mathbb{R}\nonumber\\
\gamma(t,s,x,u,z,\zeta)&=\gamma(t,s,x,u,z,\zeta,\omega):
[0,T]\times [0,T] \times \mathbb{R} \times \mathbb{U} \times \mathbb{R} \times \mathbb{R} \times \Omega \mapsto \mathbb{R}\nonumber\\
\end{align}
are given bounded $C^1$ functions with respect to $t,x$ and $u$ and $\mathbb{F}$-adapted processes in $(s,\omega)$ for each given $t,x,u,z,\zeta$.
The \emph{performance functional} $J(u)$ of a control process $u \in \mathcal{A}$ is defined by
\begin{equation}\label{eq1.4}
J(u)= \mathbb{E}[\int_0^T f(t,X(t),u(t,Z),Z)dt +g(X(T),Z)],
\end{equation}
where \begin{align}
&f(t,x,u,z): [0,T] \times\mathbb{R}\times \mathbb{U} \times\mathbb{R} \mapsto \mathbb{R}\nonumber\\
&g(x,z):\mathbb{R} \times \mathbb{R} \mapsto \mathbb{R}
\end{align}
are given bounded functions, $C^1$ with respect to $x$ and $u$, $f$ is $\mathbb{F}$-adapted for each $x,u,z$ and $g$ is $\mathcal{F}_T$-measurable for each $x,z$. The functions $f$ and $g$ are called the \emph{profit rate density} and \emph{terminal payoff density}, respectively. For completeness of the presentation we allow these functions to depend explicitly on the future value $Z$ also, although this would not be the typical case in applications. But it could be that $f$ and $g$ are influenced by the future value $Z$ directly through the action of an insider, in addition to being influenced indirectly through the control process $u$ and the corresponding state process $X$.\\

\begin{problem}
Find $u^{\star} \in\mathcal{A}$ such that
\begin{equation}\label{eq1.5}
    \sup_{u\in\mathcal{A}}J(u)=J(u^{\star}).
\end{equation}
\end{problem}

\section{The Donsker delta functional}
To study this problem we adapt the technique of the paper \cite{DO1} to the stochastic Volterra equation (SVE)  and we combine this with the method for  optimal control of SVE developed in \cite{AO}. We first recall briefly the definition and basic properties of the Donsker delta functional:

\begin{definition}
Let $Z:\Omega\rightarrow\mathbb{R}$ be a random variable which also belongs to $(\mathcal{S})^{\ast}$. Then a continuous functional
\begin{equation}\label{donsker}
    \delta_{Z}(.): \mathbb{R}\rightarrow (\mathcal{S})^{\ast}
\end{equation}
is called a Donsker delta functional of $Z$ if it has the property that
\begin{equation}\label{donsker property }
    \int_{\mathbb{R}}g(z)\delta_{Z}(z)dz= g(Z) \quad a.s.
\end{equation}
for all (measurable) $g : \mathbb{R} \rightarrow \mathbb{R}$ such that the integral converges.
\end{definition}

For example, consider the special case when $Z$ is a first order chaos random variable of the form
\begin{equation}\label{eq2.5}
    Z = Z (T_0); \text{ where } Z (t) =\int_0^t\beta(s)dB(s)+\int_0^t\int_{\mathbb{R}}\psi(s,\zeta)\tilde{N}(ds,d\zeta), \mbox{ for } t\in [0,T_0]
\end{equation}
for some deterministic functions $\beta \neq 0, \psi$ such that
\begin{equation}\label{}
    \int_0^{T_0} \{ \beta^2(t)+\int_{\mathbb{R}}\psi^2(t,\zeta)\nu(d\zeta)\} dt<\infty \text{ a.s. }
\end{equation}
and for every $\epsilon >0$ there exists $\rho > 0$ such that
\begin{equation}\label{eq8.4}
\int_{\mathbb{R} \setminus (-\epsilon,\epsilon)} e^{\rho  \zeta} \nu(d\zeta) < \infty.\nonumber\\
\end{equation}

This condition implies that the polynomials are dense in $L^2(\mu)$, where $d\mu(\zeta)=\zeta^2 d\nu(\zeta)$. It also guarantees that the measure $\nu$ integrates all polynomials of degree $\geq 2$.\\
In this case it is well known (see e.g. \cite{MOP}, \cite{DiO1}, Theorem 3.5, and \cite{DOP},\cite{DiO2}) that the Donsker delta functional exists in $(\mathcal{S})^{\ast}$ and is given
by
\begin{eqnarray}\label{eq2.7}
   \delta_Z(z)&=&\frac{1}{2\pi}\int_{\mathbb{R}}\exp^{\diamond}\big[ \int_0^{T_0}\int_{\mathbb{R}}(e^{ix\psi(s,\zeta)}-1)\tilde{N}(ds,d\zeta)+ \int_0^{T_0}ix\beta(s)dB(s)  \nonumber\\
   &+&  \int_0^{T_0}\{\int_{\mathbb{R}}(e^{ix\psi(s,\zeta)}-1-ix\psi(s,\zeta))\nu(d\zeta)-\frac{1}{2}x^2\beta^2(s)\}ds-ixz\big]dx,
\end{eqnarray}
where $\exp^{\diamond}$ denotes the Wick exponential.
Moreover, we have for $t< T_0$
\begin{align}
&\mathbb{E}[\delta_Z(z)|\mathcal{F}_t]\nonumber\\
=&  \frac{1}{2\pi}\int_{\mathbb{R}}\exp\big[\int_0^t\int_{\mathbb{R}}ix\psi(s,\zeta)\tilde{N}(ds,d\zeta) +\int_0^t ix\beta(s)dB(s)\\
&+\int_t^{T_0}\int_{\mathbb{R}}(e^{ix\psi(s,\zeta)}-1-ix\psi(s,\zeta))\nu(d\zeta)ds-\int_t^{T_0}\frac{1}{2}x^2\beta^2(s)ds-ixz\big]dx.
\end{align}

If $D_t$ and $D_{t,\zeta}$ denotes the \emph{Hida-Malliavin derivative} at $t$ and $t,\zeta$  with respect to $B$ and $\tilde{N}$, respectively, we have
\begin{align}
&\mathbb{E}[D_t\delta_Z(z)|\mathcal{F}_t]=\nonumber\\
&\frac{1}{2\pi}\int_{\mathbb{R}}\exp\big[\int_0^t\int_{\mathbb{R}}ix\psi(s,\zeta)\tilde{N}(ds,d\zeta) +\int_0^t ix\beta(s)dB(s)\nonumber\\
&+\int_t^{T_0}\int_{\mathbb{R}}(e^{ix\psi(s,\zeta)}-1-ix\psi(s,\zeta))\nu(d\zeta)ds-\int_t^{T_0}\frac{1}{2}x^2\beta^2(s)ds-ixz\big]ix\beta(t)dx
\end{align}
   and
\begin{align}
&\mathbb{E}[D_{t,z}\delta_Z(z)|\mathcal{F}_t]=\nonumber\\
& \frac{1}{2\pi}\int_{\mathbb{R}}\exp\big[\int_0^t\int_{\mathbb{R}}ix\psi(s,\zeta)\tilde{N}(ds,d\zeta) +\int_0^t ix\beta(s)dB(s)\nonumber\\
&+\int_t^{T_0}\int_{\mathbb{R}}(e^{ix\psi(s,\zeta)}-1-ix\psi(s,\zeta))\nu(d\zeta)ds-\int_t^{T_0}\frac{1}{2}x^2\beta^2(s)ds-ixz\big](e^{ix\psi(t,z)}-1)dx.
\end{align}

For more information about the Donsker delta functional, Hida-Malliavin calculus and their properties, see \cite{DO1}.

From now on we assume that $Z$ is a given random variable which also belongs to $(\mathcal{S})^{\ast}$, with a Donsker delta functional $\delta_Z(z)\in(\mathcal{S})^{\ast}$
satisfying
\begin{equation}
\mathbb{E}[\delta_Z(z)|\mathcal{F}_T] \in \mathbf{L}^2(\mathcal{F}_T,P)
\end{equation}
and
\begin{equation}
\mathbb{E}[\int_0^T(\mathbb{E}[D_t\delta_Z(z)|\mathcal{F}_t])^2 dt ]<\infty, \text{ for all } z.
\end{equation}
\section{Transforming the insider control problem to a related parameterized non-insider problem}
Since $X(t)$ is $\mathbb{H}$-adapted, we get by using the definition of the Donsker delta functional $\delta_Z(z)$ of $Z$ that
\begin{equation}\label{eq1.6}
X(t)=X(t,Z)=X(t,z)_{z=Z}=\int_{\mathbb{R}}X(t,z)\delta_Z(z)dz
\end{equation}
for some $z$-parameterized process $X(t,z)$ which is $\mathbb{F}$-adapted for each $z$.
Then, again by the definition of the Donsker delta functional we can write
\begin{align}\label{eq1.7}
&X(t)=\xi(t,Z) +\int_0^ t b (t, s, X(s,Z), u(s,Z)) ds\nonumber\\
&+ \int_0^ t \sigma (t, s, X(s,Z), u(s,Z)) dB(s) + \int_0^ t\int_{\mathbb{R}} \gamma(t, s, X(s,Z), u(s,Z), \zeta) \tilde{N}(ds, d\zeta)\nonumber\\
&=\xi(t,z)_{z=Z}+\int_0^tb (t, s, X(s,z), u(s,z))_{z=Z}ds \nonumber\\
&+ \int_0^t \sigma (t, s, X(s,z), u(s,z))_{z=Z}dB(s)\nonumber\\
&+\int_0^t \int_{\mathbb{R}} \gamma(t, s, X(s,z), u(s,z), \zeta)_{z=Z}\tilde{N}(ds,d\zeta)\nonumber\\
&= \int_{\mathbb{R}} \xi(t,z)\delta_Z(z)dz+\int_0^t \int_{\mathbb{R}}b (t, s, X(s,z), u(s,z))\delta_Z(z)dzds\nonumber\\
&+ \int_0^t \int_{\mathbb{R}}\sigma (t, s, X(s,z), u(s,z))\delta_Z(z)dzdB(s)\nonumber\\
&+\int_0^t \int_{\mathbb{R}}\int_{\mathbb{R}}\gamma(t, s, X(s,z), u(s,z), \zeta)\delta_Z(z)dz\tilde{N}(ds,d\zeta)\nonumber\\
&=\int_{\mathbb{R}} \{\xi(t,z)+\int_0^tb (t, s, X(s,z), u(s,z))ds + \int_0^t \sigma (t, s, X(s,z), u(s,z))dB(s)\nonumber\\
&+\int_0^t \int_{\mathbb{R}}\gamma(t, s, X(s,z), u(s,z), \zeta) \tilde{N}(ds,d\zeta)\} \delta_Z(z)dz.
\end{align}
Comparing \eqref{eq1.6} and \eqref{eq1.7} we see that  \eqref{eq1.6} holds if we for each $z$ choose $X(t,z)$ as the solution of the classical (but parameterized ) SVIE
\begin{align}\label{sve}
&dX(t,z) = \xi'(t,z)dt + b (t, t, X(t,z), u(t,z)) dt + (\int_0^t \frac{\partial b}{\partial t} (t, s, X(s,z), u(s,z)) ds) dt\nonumber\\
&+\sigma(t, t, X(t,z), u(t,z)) dB(t) + (\int_0^t \frac{\partial \sigma}{\partial t}(t, s, X(s,z), u(s,z)) dB(s)) dt\nonumber\\
&+\int_{\mathbb{R}} \gamma(t, t, X(t,z), u(t,z),\zeta) \tilde{N}(dt, d\zeta) + (\int_0^t\int_{\mathbb{R}} \frac{\partial \gamma}{\partial t} (t, s, X(s,z), u(s,z),\zeta) \tilde{N}(ds, d\zeta)) dt.
\end{align}
As before let $\mathcal{A}$ be the given family of admissible $\mathbb{K}-$adapted controls $u$.
Then in terms of $X(t,z)$ the performance functional $J(u)$ of a control process $u \in \mathcal{A}$ defined in \eqref{eq1.4} gets the form
\begin{align} \label{eq0.13}
  J(u) &= \mathbb{E}[\int_0^T f(t, X(t,Z),u(t,Z),Z)dt +g(X(T,Z),Z)] \nonumber\\
   &= \mathbb{E}[\int_\mathbb{R}\Big\{ \int_0^T f(t, X(t,z),u(t,z),z)\mathbb{E}[\delta_Z(z)|\mathcal{F}_t]dt \nonumber \\
   & +g(X(T,z),z)\mathbb{E}[\delta_Z(z)|\mathcal{F}_T]\Big\}dz]\nonumber\\
   &= \int_{\mathbb{R}} j(u)(z) dz,
\end{align}
where
\begin{align}\label{eq1.5}
j(u)(z)&:=  \mathbb{E}[ \int_0^T f(t, X(t,z),u(t,z),z)\mathbb{E}[\delta_Z(z)|\mathcal{F}_t]dt \nonumber \\
   & +g(X(T,z),z)\mathbb{E}[\delta_Z(z)|\mathcal{F}_T].
 \end{align}

 Thus we see that to maximize $J(u)$ it suffices to maximize $j(u)(z)$ for each value of the parameter $z \in \mathbb{R}$. Therefore Problem 1.1 is transformed into the problem
\begin{problem}
For each given $z \in \mathbb{R}$ find $u^{\star}=u^{\star}(t,z) \in\mathcal{A}$ such that
\begin{equation}\label{problem2}
    \sup_{u\in\mathcal{A}}j(u)(z)=j(u^{\star})(z).
\end{equation}
\end{problem}
\section{A sufficient-type maximum principle}
In this section we will establish a sufficient maximum principle for Problem 3.1.\\
Problem 3.1 is a stochastic control problem with a standard (albeit parameterized) stochastic Volterra equation \eqref{sve} for the state process $X(t,z)$, but with a non-standard performance functional given by \eqref{eq1.5}. We can solve this problem by a modified maximum principle approach, as follows:\\
Let $L$ and $L_{\zeta}$ be the set of all stochastic processes with parameter space $[0, T ]$ and $[0, T ] \times \mathbb{R}_0$, respectively, where $\mathbb{R}_0=\mathbb{R} \setminus\{0\}$.
Define the Hamiltonian functionals:
\begin{equation}
H_0 : [0, T] \times \mathbb{R} \times \mathbb{U} \times \mathbb{R} \times \mathbb{R} \times\mathbb{R}\rightarrow \mathbb{R}
\end{equation}
and
\begin{equation}
H_1 : [0, T] \times \mathbb{R} \times\mathbb{U} \times\mathbb{R} \times L \times L_{\zeta}\rightarrow\mathbb{R}
\end{equation}
by
\begin{align}
&H_0(t, x, z, u, p, q, r) := f(t, x, z, u)E[\delta_Z(z)|\mathcal{F}_t] + b(t, t, x, z, u)p + \sigma(t, t, x, z, u)q\nonumber\\
&+\int_{\mathbb{R}} \gamma(t, t, x, z, u,\zeta)r(\zeta)\nu(d\zeta)
\end{align}
and
\begin{align}
&H_1(t, x, z, u, p, D_tp(.), D_{t,\zeta}p(.)) := \int_t^T\frac{\partial b}{\partial s}(s, t, x, z, u)p(s,z)ds + \int_t^T\frac{\partial \sigma}{\partial s} (s, t, x, z, u)E[D_tp(s,z) | \mathcal{F}_t]ds\nonumber\\
&+ \int_t^T\int_{\mathbb{R}}\frac{\partial \gamma}{\partial s} (s, t, x, z, u, \zeta)E[D_{t,\zeta}p(s,z) | \mathcal{F}_t]\nu(d\zeta)ds.
\end{align}

Here $\mathcal{R}$ denotes the set of all functions $r(\cdot) : \mathbb{R}\rightarrow  \mathbb{R}$
such that the last integral above converges.The quantities $p,q,r(\cdot)$ are called the \emph{adjoint variables}.
Define
\begin{equation}\label{eq4.5}
H(t, x, z, u, p(.), q(.), r(.)) := H_0(t, x, z, u, p, q, r) + H_1(t, x, u, p(.), D_tp(.), D_{t,\zeta}p(.)).
\end{equation}
The \emph{adjoint processes} $p(t,z),q(t,z),r(t,z,\zeta)$ are defined as the solution of the $z$-parameterized backward stochastic differential equation (BSDE)
\begin{equation}\label{eq4.2a}
    \left\{
\begin{array}{l}
    dp(t,z) = -\frac{\partial H}{\partial x}(t,z)dt +q(t,z)dB(t)+\int_{\mathbb{R}}r(t,z,\zeta)\tilde{N}(dt,d\zeta); \\ \quad (t,z)\in(0,T)\times \mathbb{R}\\
    p(T,z)  = \frac{\partial g}{\partial x}(X(T,z),z) \mathbb{E}[\delta_Z(z)|\mathcal{F}_T];\quad z\in \mathbb{R}
\end{array}
    \right.
\end{equation}
where
\begin{equation}
\frac{\partial H}{\partial x}(t,z)=\frac{\partial H}{\partial x}(t,X(t,z),u(t,z),z,p(t,z),q(t,z),r(t,z,.)).
\end{equation}

We can now state the first maximum principle for our problem \eqref{problem2}:

\begin{theorem}{[Sufficient-type maximum principle]}\\
Let $\hat{u} \in \mathcal{A}$, and denote the associated solution  of \eqref{sve} and \eqref{eq4.2a} by $\hat{X}(t,z)$ and \\
$(\hat{p}(t,z),\hat{q}(t,z),\hat{r}(t,z,\zeta))$, respectively.  Assume that the following hold:
\begin{enumerate}
 \item $ x \rightarrow g(x,z)$ is concave for all $z$
 \item $(x,u)\rightarrow H(t,x,u,z,\widehat{p}(t,z),\widehat{q}(t,z),\hat{r}(t,z,\zeta))$ is concave for all $t,z,\zeta$
 \item $\sup_{w\in\mathbb{U}}E[H\big(t,\hat{X}(t,z),w,\widehat{p}(t,z),\widehat{q}(t,z),\hat{r}(t,z,\zeta)\big)|\mathcal{G}_t]$\\
      $=E[H\big(t,\widehat{X}(t,z),\widehat{u}(t,z),\widehat{p}(t,z),\widehat{q}(t,z),\hat{r}(t,z,\zeta)\big)|\mathcal{G}_t]$ for all $t,z,\zeta.$
 \end{enumerate}
Then $\widehat{u}(\cdot,z)$ is an optimal insider control for Problem 3.1.
\end{theorem}

\dproof  By considering an increasing sequence of stopping times $\tau_n$ converging to $T$, we may assume that all local integrals appearing in the computations below are martingales and hence have expectation 0. See \cite{OS2}. We omit the details.\\
Choose arbitrary $u(.,z)\in\mathcal{A}$, and let the corresponding solution of \eqref{sve} and \eqref{eq4.2a}  be $X(t,z)$, $p(t,z)$, $q(t,z)$, $r(t,z,\zeta)$.
For simplicity of notation we write\\
$f=f(t,X(t,z),u(t,z),z)$,
$\widehat{f}=f(t,\widehat{X}(t,z),\widehat{u}(t,z),z)$
 and similarly with $b$, $\widehat{b}$, $\sigma$, $\widehat{\sigma}$ and so on.\\
 Moreover put
 \begin{equation}
 \hat{H}(t)=H(t,\widehat{X}(t,z),\widehat{u}(t,z),\widehat{p}(t,z),\widehat{q}(t,z),\widehat{r}(t,z,.))
 \end{equation}
 and
 \begin{equation}
 H(t)=H(t,X(t,z), u(t,z),\widehat{p}(t,z),\widehat{q}(t,z),\widehat{r}(t,z,.))
 \end{equation}
 and similarly with $\hat{H}_0(t)$ and  $H_0(t)$.
  In the following we write  $\widetilde{f}=f-\widehat{f}$, $\widetilde{b}=b-\widehat{b}$, $\widetilde{X}=X-\widehat{X}$.\\
 Consider
 \begin{equation*}
    j(u(.,z))-j(\widehat{u}(.,z))=I_1+ I_2,
 \end{equation*}
 where
 \begin{equation}\label{eq4.7}
    I_1=\mathbb{E}[\int_0^T\{f(t)-\widehat{f}(t)\}\mathbb{E}[\delta_Z(z)|\mathcal{F}_t]dt], \quad I_2=\mathbb{E}[\{g(x)-\hat{g}(x)\}\mathbb{E}[\delta_Z(z)|\mathcal{F}_T]].
 \end{equation}
  By the definition of $H$ we have
  \begin{align}\label{eq4.8}
    I_1 &= \mathbb{E}[\int_0^T\{H_0(t)-\widehat{H_0}(t)-\widehat{p}(t)\widetilde{b}(t) - \widehat{q}(t)\widetilde{\sigma}(t)\nonumber\\
   & -\int_{\mathbb{R}}\hat{r}(t,\zeta)\tilde{\gamma}(t,\zeta)\nu(d\zeta)\} dt].
  \end{align}
   Since $g$ is concave with respect to $x$ we have
   \begin{align}
&(g(X(T,z),z)-g(\hat{X}(T,z),z))\mathbb{E}[\delta_Z(z)|\mathcal{F}_T]\nonumber\\
&\leq \frac{\partial g}{\partial x}(\hat{X}(T,z),z)\mathbb{E}[\delta_Z(z)|\mathcal{F}_T](X(T,z)-\hat{X}(T,z)),
\end{align}
and hence
 \begin{align}\label{eq4.11}
I_2& \leq\mathbb{E}\Big[\frac{\partial g}{\partial x}(\widehat{X}(T,z))\mathbb{E}[\delta_Z(z)|\mathcal{F}_T]\tilde{X}(T,z)\Big]=\mathbb{E}[\widehat{p}(T,z)\widetilde{X}(T,z)] \\ \nonumber
  & = E \Big[\int_0^T \Big\{\hat{p}(t,z) \Big(b(t,z) -\hat{b}(t,z) + \int_0^t(\frac{\partial b}{\partial t} (t, s)-\frac{\partial \hat{b}}{\partial t} (t, s)) ds\\ \nonumber
&+ \int_0^t (\frac{\partial \sigma}{\partial t} (t, s) - \frac{\partial \hat{\sigma}}{\partial t} (t, s) )dB(s)\\ \nonumber
&+ \int_0^t \int_{\mathbb{R}}(\frac{\partial \gamma}{\partial t}(t, s,\zeta) -\frac{\partial\hat{\gamma}}{\partial t}(t, s,\zeta))\tilde{N}(ds, d\zeta)\Big)\\ \nonumber
&-\frac{\partial \hat{H}}{\partial x}(t)(X(t,z)-\hat{X}(t,z))\\ \nonumber
&+\hat{q}(t,z)[\sigma (t,z) - \hat{\sigma}(t,z)] + \int_{\mathbb{R}} \hat{r}(t,z,\zeta)[\gamma(t,z,\zeta)-\hat{\gamma}(t,z,\zeta)]\nu(d\zeta)\Big\}dt\Big]
\end{align}
By the Fubini Theorem, we get
\begin{equation}\label{eq4.14vol}
\int_0^T (\int_0^t \frac{\partial b}{\partial t}(t, s)ds)\hat{p}(t,z)dt = \int_0^T (\int_s^T\frac{\partial b}{\partial t}(t, s)\hat{p}(t,y)dt)ds =\int_0^T (\int_t^T \frac{\partial b}{\partial s}(s, t)\hat{p}(s,z)ds)dt,
\end{equation}
and similarly, by the duality theorems,
\begin{align}\label{eq4.15vol}
&E[\int_0^T (\int_0^t \frac{\partial \sigma}{\partial t}(t, s)dB(s))\hat{p}(t,z)dt] = \int_0^T E[\int_0^t\frac{\partial\sigma}{\partial t}(t, s)dB(s) \hat{p}(t,z)]dt\nonumber\\
&=\int_0^T E[\int_0^t \frac{\partial \sigma}{\partial t}(t, s)E[D_s\hat{p}(t,z)|\mathcal{F}_s]ds]dt
= \int_0^T E[\int_s^T \frac{\partial \sigma}{\partial t}(t, s)E[D_s\hat{p}(t,z)|\mathcal{F}_s]dt]ds \nonumber\\
&= E[\int_0^T \int_t^ T \frac{\partial \sigma}{\partial s}(s, t)E[D_t\hat{p}(s,z) |\mathcal{F}_t]dsdt]
\end{align}

and
\begin{align}\label{eq4.16vol}
&E[\int_0^T (\int_0^t \int_{\mathbb{R}} \frac{\partial \gamma}{\partial t}(t, s,\zeta)\tilde{N}(ds, d\zeta))\hat{p}(t,z)dt] = \int_0^T E[\int_0^t\int{\mathbb{R}}\frac{\partial \gamma}{\partial t} (t, s, \zeta)\tilde{N}(ds, d\zeta))\hat{p}(t,z)]dt\nonumber\\
&=\int_0^T E[\int_0^t\int_{\mathbb{R}} \frac{\partial \gamma}{\partial t} (t, s,\zeta)E[D_{s,\zeta}\hat{p}(t,z)|\mathcal{F}_s]\nu(d\zeta)ds]dt\nonumber\\
&= \int_0^T E[\int_s^T \int_{\mathbb{R}} \frac{\partial \gamma}{\partial t}(t, s,\zeta)E[D_{s,\zeta}\hat{p}(t,z) | \mathcal{F}_s]\nu(d\zeta)dt]ds\nonumber\\
&= E[\int_0^T\int_t^T\int_{\mathbb{R}} \frac{\partial \gamma}{\partial s} (s, t,\zeta)E[D_{t,\zeta}\hat{p}(s,z)|\mathcal{F}_t]\nu(d\zeta)dsdt].
\end{align}

Substituting \eqref{eq4.14vol}, \eqref{eq4.15vol} and \eqref{eq4.16vol} into \eqref{eq4.11}, we get
\begin{align}\label{eq4.17}
&I_2 \leq E \Big[\int_0^T \Big(\hat{p}(t,z)[b(t) - \hat{b}(t)] +\hat{p}(t,z) \int_t^T(\frac{\partial b}{\partial s} (s,t) -\frac{\partial \hat{b}}{\partial s} (s, t))ds\nonumber\\
&+\int_t^T (\frac{\partial \sigma}{\partial s} (s, t) - \frac{\partial \hat{\sigma}}{\partial s}(s, t))E[D_t\hat{p}(s,z) |\mathcal{F}_t]ds\nonumber\\
&+ \int_t^T\int_{\mathbb{R}}(\frac{\partial \gamma}{\partial s}(s, t,\zeta) -\frac{\partial \hat{\gamma}}{\partial s}(s, t,\zeta))E[D_{t,\zeta}\hat{p}(s)| \mathcal{F}_t]ds\nu(d\zeta)\nonumber\\
&-\frac{\partial\hat{H}}{\partial x}(t)(X(t,z) - \hat{X}(t,z))\nonumber\\
&+\hat{q}(t,z) [\sigma(t,z) - \hat{\sigma}(t,z)] + \int_{\mathbb{R}} \hat{r}(t,\zeta) [\gamma(t,z,\zeta) - \hat{\gamma}(t,z,\zeta)]\nu(d\zeta)\Big)dt\Big] .
\end{align}
Adding \eqref{eq4.8} and \eqref{eq4.17}, we get
\begin{align}\label{eq4.10}
&J(u) - J(\hat{u}) = I_1 + I_2\nonumber\\
&\leq E \Big[\int_0^T \Big\{H(t) - \hat{H}(t) -\frac{\partial{H}}{\partial x}(t)(X(t,z) - \hat{X}(t,z))\Big\} dt\Big].
\end{align}

By the concavity assumption of $H$ in $(x,u)$ we have:
\begin{equation}
H(t)-\hat{H}(t)\leq \frac{\partial\widehat{H}}{\partial x}(t)(X-\hat{X})(t,z)+ \frac{\partial \widehat{H}}{\partial u}(t)(u(t,z)-\hat{u}(t,z)),
\end{equation}
Then equation \eqref{eq4.10} becomes
\begin{align}\label{eq4.10a}
J(u) - J(\hat{u}) &\leq E \Big[\int_0^T\frac{\partial \widehat{H}}{\partial u}(t)(u(t,z)-\hat{u}(t,z)) dt\Big]\nonumber\\
&=\int_0^TE \Big[\frac{\partial \widehat{H}}{\partial u}(t)(u(t,z)-\hat{u}(t,z)) \Big]dt\nonumber\\
&=\int_0^TE \Big[E[\frac{\partial \widehat{H}}{\partial u}(t)(u(t,z)-\hat{u}(t,z))|\mathcal{G}_t] \Big]dt
\end{align}

and the maximum condition implies that
\begin{equation}
E[\frac{\partial\widehat{H}}{\partial u}(t,z)(u(t,z)-\hat{u}(t,z))|\mathcal{G}_t]\leq 0.
\end{equation}
Hence by \eqref{eq4.10a} we get
$j(u)\leq j(\hat{u})$.
Since $u\in\mathcal{A}$ was arbitrary, this shows that $\hat{u}$ is optimal.
\fproof

\section{A necessary-type maximum principle}

We proceed to establish a corresponding necessary maximum principle. For this, we do not need concavity conditions, but instead we need the following assumptions about the set of admissible control processes:\\
\begin{itemize}
\item
$A_1$. For all $t_0\in [0,T]$ and all bounded $\mathcal{G}_{t_0}$-measurable random variables $\alpha(z,\omega)$, the control
$\theta(t,z, \omega) := \mathbf{1}_{[t_0,T ]}(t)\alpha(z,\omega)$ belongs to $\mathcal{A}$.\\
\item
$A_2$. For all $u, \beta_0 \in\mathcal{A}$ with $\beta_0(t,z) \leq K < \infty$ for all $t,z$  define
\begin{equation}\label{delta}
    \delta(t,z)=\frac{1}{2K}dist(u(t,z),\partial\mathbb{U})\wedge1 > 0
\end{equation}
and put
\begin{equation}\label{eq3.2}
    \beta(t,z)=\delta(t,z)\beta_0(t,z).
\end{equation}
Then the control
\begin{equation*}
    \widetilde{u}(t,x,z)=u(t,z) + a\beta(t,z) ; \quad t \in [0,T]
\end{equation*}
belongs to $\mathcal{A}$ for all $a \in (-1, 1)$.\\
\item
$A3$. For all $\beta$ as in \eqref{eq3.2} the derivative process
\begin{equation}\label{eq5.3a}
    \chi(t,z):=\frac{d}{da}X^{u+a\beta}(t,z)|_{a=0}
\end{equation}
exists, and belong to $\mathbf{L}^2(\lambda\times \mathbf{P})$ and
\begin{align}
&\chi(t,z) =\int^t_0\Big(\frac{\partial b}{\partial x}(t, s)\chi(s,z) + \frac{\partial b}{\partial u}(t, s)\beta(s,z)\Big)ds\nonumber\\
&+\int^t_0\Big(\frac{\partial \sigma}{\partial x}(t, s)\chi(s,z) + \frac{\partial \sigma}{\partial u}(t, s)\beta(s,z)\Big)dB(s)\nonumber\\
&+\int^t_0\int_{\mathbb{R}}\Big(\frac{\partial \gamma}{\partial x}(t, s, \zeta)\chi(s,z) + \frac{\partial \gamma}{\partial u}(t, s, \zeta)\beta(s,z)\Big)  \tilde{N}(ds, d\zeta),
\end{align}

and hence
\begin{align}\label{d chi}
&d\chi(t,z) =\Big[\frac{\partial b}{\partial x}(t, t)\chi(t,z) +\frac{\partial b}{\partial u} (t, t)\beta(t,z) +
\int_0^t\Big(\frac{\partial^2 b}{\partial t\partial x} (t, s)\chi(s,z) + \frac{\partial^2 b}{\partial t\partial u} (t, s)\beta(s,z)\Big)ds\nonumber\\
&+\int_0^t\Big(\frac{\partial^2 \sigma}{\partial t\partial x} (t, s)\chi(s,z) + \frac{\partial^2 \sigma}{\partial t\partial u} (t, s)\beta(s,z)\Big)dB(s)\nonumber\\
&+\int_0^t\int_{\mathbb{R}}\Big(\frac{\partial^2 \gamma}{\partial t\partial x} (t, s,\zeta)\chi(s,z) + \frac{\partial^2 \gamma}{\partial t\partial u} (t, s,\zeta)\beta(s,z)\Big) \tilde{N}(ds, d\zeta)\Big]dt\nonumber\\
&+\Big(\frac{\partial \sigma}{\partial x}(t, t)\chi(t,z) +\frac{\partial \sigma}{\partial u} (t, t)\beta(t,z)\Big) dB(t)\nonumber\\
&+\int_{\mathbb{R}}\Big(\frac{\partial \gamma}{\partial x}(t, t)\chi(t,z) +\frac{\partial \gamma}{\partial u} (t, t)\beta(t,z)\Big)\tilde{N}(dt, d\zeta)
\end{align}
\end{itemize}

\begin{theorem}{[Necessary-type maximum principle]}\label{necessary theorem} \\
Let $\hat{u} \in \mathcal{A}$ and $z \in \mathbb{R}$. Then the following are equivalent:
\begin{enumerate}
\item $\frac{d}{da}j(\hat{u}+a\beta)(z)|_{a=0}=0$ for all bounded $\beta \in \mathcal{A}$ of the form \eqref{eq3.2}.
\item $E[\frac{\partial H}{\partial u}(t,z)|\mathcal{G}_t]_{u=\hat{u}}=0$ for all $t\in [0,T].$
\end{enumerate}
\end{theorem}

\dproof
For simplicity of notation we write $u$ instead of $\hat{u}$ in the following. \\
By considering an increasing sequence of stopping times $\tau_n$ converging to $T$, we may assume that all local integrals appearing in the computations below are martingales and have expectation 0. See \cite{OS2}. We omit the details.\\
We can write
$$\frac{d}{da}j((u+a\beta)(z))|_{a=0}=I_1+I_2$$
where
$$I_1=\frac{d}{da}\mathbb{E}[\int_0^Tf(t,X^{u+a\beta}(t,z),u(t,z)+a\beta(t,z),z)\mathbb{E}[\delta_Z(z)|\mathcal{F}_t]dt]|_{a=0}$$
and
$$I_2=\frac{d}{da}\mathbb{E}[g(X^{u+a\beta}(T,z),z)\mathbb{E}[\delta_Z(z)|\mathcal{F}_T]]|_{a=0}.$$
By our assumptions on $f$ and $g$ and by \eqref{eq5.3a} we have
\begin{equation}\label{iii0}
    I_1=\mathbb{E}[\int_0^T\{\frac{\partial f}{\partial x}(t,z)\chi(t,z)+\frac{\partial f}{\partial u}(t,z)\beta(t,z)\}\mathbb{E}[\delta_Z(z)|\mathcal{F}_t]dt],
\end{equation}
\begin{equation}\label{iii2}
    I_2=\mathbb{E}[\frac{\partial g}{\partial x}(X(T,z),z)\chi(T,z)\mathbb{E}[\delta_Z(z)|\mathcal{F}_T]]=\mathbb{E}[p(T,z)\chi(T,z)].
\end{equation}
We have
\begin{align}
&I_1=\mathbb{E}[\int_0^T\Big(\{ \frac{\partial H_0}{\partial x}(t)-\frac{\partial b}{\partial x}(t,t)p(t,z)-\frac{\partial \sigma}{\partial x}(t,t)q(t,z)\nonumber\\
&-\int_{\mathbb{R}}\frac{\partial \gamma}{\partial x}(t,t,\zeta)r(t,z,\zeta)\nu(d\zeta)\}X(t,z)+\{ \frac{\partial H_0}{\partial u}-\frac{\partial b}{\partial u}(t,t)p(t,z)-\frac{\partial \sigma}{\partial u}(t,t)q(t,z)\nonumber\\
&-\int_{\mathbb{R}}\frac{\partial \gamma}{\partial u}(t,t,\zeta)r(t,z,\zeta)\nu(d\zeta)\}\beta(t,z)\Big)dt]
\end{align}
Since $H_0=H-H_1$ then
\begin{align}\label{iii1}
&I_1=\mathbb{E}[\int_0^T\Big(\{ \frac{\partial H}{\partial x}(t)-\frac{\partial H_1}{\partial x}(t)-\frac{\partial b}{\partial x}(t,t)p(t,z)-\frac{\partial \sigma}{\partial x}(t,t)q(t,z)\nonumber\\
&-\int_{\mathbb{R}}\frac{\partial \gamma}{\partial x}(t,t,\zeta)r(t,z,\zeta)\nu(d\zeta)\}\chi(t,z)+\{ \frac{\partial H}{\partial u}(t)-\frac{\partial H_1}{\partial u}(t)-\frac{\partial b}{\partial u}(t,t)p(t,z)-\frac{\partial \sigma}{\partial u}(t,t)q(t,z)\nonumber\\
&-\int_{\mathbb{R}}\frac{\partial \gamma}{\partial u}(t,t,\zeta)r(t,z,\zeta)\nu(d\zeta)\}\beta(t,z)\Big)dt]\nonumber\\
&=\mathbb{E}[\int_0^T\Big(\{\frac{\partial H}{\partial x}-\int_t^T\frac{\partial^2b}{\partial s\partial x}(s,t)p(s,z)ds-\int_t^T\frac{\partial^2\sigma}{\partial s\partial x}(s,t)E[D_tp(s,z)|\mathcal{F}_t]ds\nonumber\\
&-\int_t^T\int_{\mathbb{R}}\frac{\partial^2\gamma}{\partial s\partial x}(s,t)E[D_{t,\zeta}p(s,z)|\mathcal{F}_t]\nu(d\zeta)ds-\frac{\partial b}{\partial x}(t,t)p(t,z)-\frac{\partial \sigma}{\partial x}(t,t)q(t,z)\nonumber\\
&-\int_{\mathbb{R}}\frac{\partial \gamma}{\partial x}(t,t,\zeta)r(t,z,\zeta)\nu(d\zeta)\}\chi(t,z)\nonumber\\
&+\{\frac{\partial H}{\partial u}(t)-\int_t^T\frac{\partial^2b}{\partial s\partial u}(s,t)p(s,z)ds-\int_t^T\frac{\partial^2\sigma}{\partial s\partial u}(s,t)E[D_tp(s,z)|\mathcal{F}_t]ds\nonumber\\
&-\int_t^T\int_{\mathbb{R}}\frac{\partial^2\gamma}{\partial s\partial u}(s,t)E[D_{t,\zeta}p(s,z)|\mathcal{F}_t]\nu(d\zeta)ds-\frac{\partial b}{\partial u}(t,t)p(t,z)\nonumber\\
&-\frac{\partial \sigma}{\partial u}(t,t)q(t,z)
-\int_{\mathbb{R}}\frac{\partial \gamma}{\partial u}(t,t,\zeta)r(t,z,\zeta)\nu(d\zeta)\}\beta(t,z)\Big)dt]
\end{align}

By the It\^{o} formula

\begin{align}\label{eq5.9a}
   I_2&= \mathbb{E}[p(T,z)\chi(T,z)]=\mathbb{E}[\int_0^Tp(t,z)d\chi(t,z) +\int_0^T\chi(t,z)dp(t,z)\nonumber\\
   &+\int_0^Td[\chi,p](t,z)] \nonumber\\
      &= E\Big[\int_0^Tp(t)(\frac{\partial b}{\partial x} (t, t)\chi(t,z) + \frac{\partial b}{\partial u}(t, t)\beta(t,z)) dt\nonumber\\
&+\int^T_0p(t)\Big\{\int^t_0(\frac{\partial^2 b}{\partial t \partial x}(t, s)\chi(s,z) +\frac{\partial^2 b}{\partial t \partial u}(t, s)\beta(s,z))ds\Big\}dt\nonumber\\
&+\int^T_0p(t)\Big\{\int^t_0 (\frac{\partial^2 \sigma}{\partial t \partial x}(t, s)\chi(s,z) +\frac{\partial^2 \sigma}{\partial t \partial u}(t, s)\beta(s,z))dB(s)\Big\}dt\nonumber\\
&+\int^T_0p(t)\Big\{\int^t_0 (\frac{\partial^2 \gamma}{\partial t \partial x}(t, s,\zeta)\chi(s,z) +\frac{\partial^2 \gamma}{\partial t \partial x}(t, s, \zeta)\beta(s,z))\tilde{N}(ds,d\zeta)\Big\}dt\nonumber\\
&-\int^T_0\chi (t,z) \frac{\partial H}{\partial x}(t)dt +\int_0^Tq(t,z)(\frac{\partial \sigma}{\partial x}(t, t)\chi(t,z) + \frac{\partial \sigma}{\partial u}(t, t)\beta(t,z))dt\nonumber\\
&+\int^T_0\int_{\mathbb{R}}r(t,\zeta)(\frac{\partial \gamma}{\partial x}(t, t,\zeta)\chi(t,z) +\frac{\partial \gamma}{\partial u}(t, t,\zeta)\beta(t,z))\nu(d\zeta)dt\Big].
\end{align}
 From \eqref{eq4.14vol}, \eqref{eq4.15vol} and \eqref{eq4.16vol}, we have
 \begin{align}
&E [p(T )Y (T )]\nonumber\\
&= E\Big[\int_0^T \Big\{\frac{\partial b}{\partial x}(t, t)p(t,z) + \int_t^T\Big(\frac{\partial^2 b}{\partial s\partial x} (s, t)p(s,z) +\frac{\partial^2 \sigma}{\partial s\partial x}(s, t)E[D_tp(s,z) | \mathcal{F}_t]\nonumber\\
&+\int_{\mathbb{R}}\frac{\partial^2 \gamma}{\partial s\partial x} (s, t, \zeta) E[D_{t,\zeta}p(s)|\mathcal{F}_t]\nu(d\zeta)\Big)ds\Big\}\chi(t,z)dt\nonumber\\
&+\int_0^T\Big\{\frac{\partial b}{\partial u}(t, t)p(t,z) + \int_t^T\Big(\frac{\partial^2 b}{\partial s\partial u} (s, t)p(s,z) +\frac{\partial^2 \sigma}{\partial s\partial u}(s, t)E[D_tp(s,z) | \mathcal{F}_t]\nonumber\\
&+\int_{\mathbb{R}}\frac{\partial^2 \gamma}{\partial s\partial u} (s, t, \zeta) E[D_{t,\zeta}p(s)|\mathcal{F}_t]\nu(d\zeta)\Big)ds\Big\}\beta(t,z)dt\nonumber\\
&-\int_0^T\frac{\partial H}{\partial x}(t)X(t,z)dt +\int_0^T(\frac{\partial \sigma}{\partial x}(t, t)\chi(t,z) + \frac{\partial \sigma}{\partial u}(t, t)\beta(t)) q(t,z)dt\nonumber\\
&+\int_0^T\int_{\mathbb{R}}(\frac{\partial \gamma}{\partial x}(t, t,\zeta)\chi(t,z) + \frac{\partial \sigma}{\partial u}(t, t,\zeta)\beta(t,z)) r(t,\zeta)\nu(d\zeta)dt\Big].
\end{align}
Summing (\ref{iii1}) and \eqref{eq5.9a} we get
\begin{equation*}
    \frac{d}{da}j((u+a\beta)(.,z))|_{a=0}=I_1+I_2=\mathbb{E}[\int_0^T\frac{\partial H}{\partial u}(t,z)\beta(t,z)dt].
\end{equation*}
We conclude that
\begin{equation*}
    \frac{d}{da}j(u+a\beta)(z))|_{a=0}=0
\end{equation*}
if and only if
\begin{equation}
\mathbb{E}[\int_0^T\frac{\partial H}{\partial u}(t,z)\beta(t,z)dt]=0,
\end{equation}
for all bounded $\beta\in\mathcal{A}$ of the form \eqref{eq3.2}.\\

\noindent In particular, applying this to $\beta(t,z) = \theta(t,z)$ as in $A1$,  we get that this is again equivalent to
\begin{equation}\label{pass}
\mathbb{E}[\alpha(z,\omega)\int_{t_0}^T\frac{\partial H}{\partial u}(t,z)dt]=0,
\end{equation}
Differentiating the right-hand side of \eqref{pass}, we get
\begin{equation}
\mathbb{E}[\alpha(z,\omega)\frac{\partial H}{\partial u}(t_0,z)]=0
\end{equation}
Since this holds for all bounded $\mathcal{G}_{t_0}$-measurable $\alpha$, so we deduce that
\begin{equation}
\mathbb{E}[\frac{\partial H}{\partial u}(t_0,z)|\mathcal{F}_{t_0}]=0, \forall t_0\in[0,T]
\end{equation}
\fproof
\section{Applications}
\subsection{The case when the coefficients do not depend on $x$}
Consider the case when the coefficients do not depend on $x$, i.e., the system has the form:
\begin{align}\label{eq6.1}
X(t,z) = \xi(t,z) + \int_0^tb(t, s, u(s,z),z)ds + \int_0^t\sigma(t, s, u(s,z),z)dB(s) + \int_0^t\int_{\mathbb{R}}\gamma(t, s, u(s,z),\zeta)\tilde{N}(ds, d\zeta)
\end{align}
with performance functional
\begin{equation}\label{eq6.2}
J(u) = E[ \int_0^T f(t, u(t))E[\delta_Z(z)|\mathcal{F}_t]dt + g(X(T))E[\delta_Z(z)|\mathcal{F}_T]] .
\end{equation}

In this case the Hamiltonian $H$ given in \eqref{eq4.5} takes the form
\begin{align}\label{eq6.3}
&H(t, z, u, p, q, r) := f(t,z, u)E[\delta_Z(z)|\mathcal{F}_t] + b(t, t, z, u)p + \sigma(t, t,z, u)q +\int_{\mathbb{R}} \gamma(t, t,z, u,\zeta)r(\zeta)\nu(d\zeta)\nonumber\\
&\int_t^T\frac{\partial b}{\partial s}(s, t,z, u)p(s)ds + \int_t^T\frac{\partial \sigma}{\partial s} (s, t, z, u)E[D_tp(s) | \mathcal{F}_t]ds
+ \int_t^T\int_{\mathbb{R}}\frac{\partial \gamma}{\partial s} (s, t, z, u, \zeta)E[D_{t,\zeta}p(s) | \mathcal{F}_t]\nu(d\zeta)ds
\end{align}
The BSDE \eqref{eq4.2a} for the adjoint variables $p, q, r$ gets the form
\begin{equation}
\begin{cases}
dp(t,z) = q(t,z)dB(t) + \int_{\mathbb{R}}r(t,z,\zeta)\tilde{N}(dt, d\zeta), 0 \leq t \leq T\\
p(T,z) = g'(X(T))E[\delta_Z(z)|\mathcal{F}_T]
\end{cases}
\end{equation}
which has the solution
\begin{align}
&p(t,z) = E[g'(X(T,z)) E[\delta_Z(z)|\mathcal{F}_T]|\mathcal{F}_t] \label{eq6.5}\\
&q(t,z) = D_tp(t,z) = E[D_t(g'(X(T,z))E[\delta_Z(z)|\mathcal{F}_T]) | \mathcal{F}_t] \label{eq6.6}\\
&r(t,z,\zeta) = D_{t,\zeta}p(t,z) = E[D_{t,\zeta}(g'(X(T,z))E[\delta_Z(z)|\mathcal{F}_T])|\mathcal{F}_t].\label{eq6.7}
\end{align}
Substituting \eqref{eq6.5}-\eqref{eq6.7} into \eqref{eq6.3} we get
\begin{equation}
E[H(t, z, u, p(.), q(.), r(.)) | \mathcal{F}_t] = E[H_0(t, z, u, p, q, r) | \mathcal{F}_t],
\end{equation}
where
\begin{align}
&H_0(t, z, u, p, q, r) := f(t,z, u)E[\delta_Z(z)|\mathcal{F}_t] + b(t, t, z, u) E[g'(X(T,z)) E[\delta_Z(z)|\mathcal{F}_T]|\mathcal{F}_t]\nonumber\\
&+ \sigma(t, t,z, u)E[D_t(g'(X(T,z))E[\delta_Z(z)|\mathcal{F}_T]) | \mathcal{F}_t]\nonumber\\
 &+\int_{\mathbb{R}} \gamma(t, t,z, u,\zeta)E[D_{t,\zeta}(g'(X(T,z))E[\delta_Z(z)|\mathcal{F}_T])|\mathcal{F}_t]\nu(d\zeta)\nonumber\\
&+\int_t^T\frac{\partial b}{\partial s}(s, t,z, u) g'(X(T,z)) E[\delta_Z(z)|\mathcal{F}_T]  ds\nonumber\\
 &+ \int_t^T\frac{\partial \sigma}{\partial s} (s, t, z, u)E[D_t\{g'(X(T,z)) E[\delta_Z(z)|\mathcal{F}_T]\}|\mathcal{F}_t]ds\nonumber\\
&+ \int_t^T\int_{\mathbb{R}}\frac{\partial \gamma}{\partial s} (s, t, z, u, \zeta)E[D_{t,\zeta}\{g'(X(T,z)) E[\delta_Z(z)|\mathcal{F}_T]\}|\mathcal{F}_t]\nu(d\zeta)ds
\end{align}
Since we have
\begin{align}
&\int_t^T\frac{\partial b}{\partial s}(s, t,z, u) g'(X(T,z)) E[\delta_Z(z)|\mathcal{F}_T]  ds\nonumber\\
&=\int_t^T\frac{\partial b}{\partial s}(s, t,z, u)ds g'(X(T,z)) E[\delta_Z(z)|\mathcal{F}_T]\nonumber\\
&=(b(T, t,z, u)-b(t, t,z, u))g'(X(T,z)) E[\delta_Z(z)|\mathcal{F}_T],
\end{align}

\begin{align}
&\int_t^T\frac{\partial \sigma}{\partial s} (s, t, z, u)E[D_t\{g'(X(T,z)) E[\delta_Z(z)|\mathcal{F}_T]\}|\mathcal{F}_t]ds\nonumber\\
&=\int_t^T\frac{\partial \sigma}{\partial s} (s, t, z, u)dsE[D_t\{g'(X(T,z)) E[\delta_Z(z)|\mathcal{F}_T]\}|\mathcal{F}_t]\nonumber\\
&=(\sigma(T, t, z, u)-\sigma(T, t, z, u))E[D_t\{g'(X(T,z)) E[\delta_Z(z)|\mathcal{F}_T]\}
\end{align}
and
\begin{align}
&\int_t^T\int_{\mathbb{R}}\frac{\partial \gamma}{\partial s} (s, t, z, u, \zeta)E[D_{t,\zeta}\{g'(X(T,z)) E[\delta_Z(z)|\mathcal{F}_T]\}|\mathcal{F}_t]\nu(d\zeta)ds\nonumber\\
&=\int_{\mathbb{R}}(\int_t^T \frac{\partial \gamma}{\partial s} (s, t, z, u, \zeta)ds)E[D_{t,\zeta}\{g'(X(T,z)) E[\delta_Z(z)|\mathcal{F}_T]\}|\mathcal{F}_t]\nu(d\zeta)\nonumber\\
&=\int_{\mathbb{R}}(\gamma(T, t, z, u, \zeta)-\gamma(t, t, z, u, \zeta))E[D_{t,\zeta}\{g'(X(T,z)) E[\delta_Z(z)|\mathcal{F}_T]\}|\mathcal{F}_t]\nu(d\zeta).
\end{align}

Then $H_0(t, v, p, q, r)$ is reduced to
\begin{align}
&\mathbb{H}_0(t, z,u, X(T,z))  =f(t,z, u)E[\delta_Z(z)|\mathcal{F}_t] + b(T, t, z, u) E[g'(X(T,z)) E[\delta_Z(z)|\mathcal{F}_T]|\mathcal{F}_t]\nonumber\\
&+ \sigma(T, t,z, u)E[D_t(g'(X(T,z))E[\delta_Z(z)|\mathcal{F}_T]) | \mathcal{F}_t]\nonumber\\
 &+\int_{\mathbb{R}} \gamma(T, t,z, u,\zeta)E[D_{t,\zeta}(g'(X(T,z))E[\delta_Z(z)|\mathcal{F}_T])|\mathcal{F}_t]\nu(d\zeta)
\end{align}
We conclude that, in this case, we have the following maximum principles:
\begin{theorem}(Sufficient maximum principle)\\
 Suppose that the coefficients $f(t, z, u), b(t, s, z, u), \sigma(t, s, z, u)$
and $\gamma(t, s, z, u, \zeta)$ of the stochastic control system \eqref{eq6.1}-\eqref{eq6.2} do not depend on $x$.
Let $\hat{u}\in\mathcal{A}$ with associated solution $\hat{X}$ of \eqref{eq6.1}. Suppose that the functions
\begin{equation}
x \rightarrow g(x)
\end{equation}
and
\begin{equation}
u\rightarrow \mathbb{H}_0(t, z, u, \hat{X}(T,z))
\end{equation}
are concave and that, for all $t$,
\begin{equation}
\max_{u\in\mathbb{U}}E[\mathbb{H}_0(t, z, u, \hat{X}(T,z))|\mathcal{G}_t]= E[\mathbb{H}_0(t, \hat{u}(t,z), \hat{X}(T,z))|\mathcal{G}_t].
\end{equation}
Then, $\hat{u}$ is an optimal control, i.e.,
\begin{equation}
\sup_{u\in A}J(u) = J(\hat{u}).
\end{equation}
\end{theorem}

\begin{theorem}(Necessary maximum principle) Let $X(t,z)$ and $J(u)$ be as in Theorem 6.1.
Let $\hat{u}\in\mathcal{A}$ with associated solution $\hat{X}$ of \eqref{eq6.1}.
Then the following, (i) and (ii), are equivalent:\\
(i)$\frac{d}{da}J(\hat{u} + a\beta)|_{a=0}= 0$
for all processes $\beta\in\mathcal{A}$ of the form \eqref{eq3.2}.\\
(ii)$E[\frac{\partial\mathbb{H}_0}{\partial u} (t, u, \hat{X}(T,z))|\mathcal{G}_t]|_{u=\hat{u}(t,z)}=0$
\end{theorem}

\section{Example: Optimal insider portfolio in a financial market modeled by a Volterra equation}
In this example, we choose $\mathbb{G} = \mathbb{F}$.
Consider a financial market where the unit price $S_0(t)$ of the risk free asset is
\begin{equation}
S_0(t) = 1, t\in [0, T]
\end{equation}
and the unit price process $S(t)$ of the risky asset has no jumps and is given by
\begin{equation}
\begin{cases}
dS(t) = b_0(t,t, Z )S(t)dt + \sigma_0(t,t, Z )S(t)dB(t)\\
+(\int_0^t\frac{\partial b_0}{\partial t}(t,s)S(s)ds+\int_0^t\frac{\partial \sigma_0}{\partial t}(t,s)S(s)dB(s))dt; t \in[0,T]\\
S(0)>0
\end{cases}
\end{equation}
Then the wealth process $X(t) = X^{\Pi}(t)$ associated to a portfolio $u(t) =\Pi(t)$, interpreted
as the fraction of the wealth invested in the risky asset at time $t$, is described by the linear stochastic Volterra equation
\begin{equation}
X(t) = x_0 + \int_0^t b_0(t, s, Z)\Pi(s)X(s)ds + \int_0^t \sigma_0(t, s,Z)\Pi(s)X(s)dB(s); t\geq0
\end{equation}
or, in differential form
\begin{equation}
\begin{cases}
dX(t) =b_0(t, t, Z)\Pi(t)X(t)dt +\sigma_0(t, t, Z)\Pi(t)X(t)dB(t)\\
+(\int_0^t\frac{\partial b_0}{\partial t}(t,s)\Pi(s)X(s)ds+\int_0^t\frac{\partial \sigma_0}{\partial t}(t,s)\Pi(s)X(s)dB(s))dt; t\in[0, T] \\
X(0)=x_0>0
\end{cases}
\end{equation}
Let $U$ be a given utility function. We want to find $\Pi\in\mathcal{A}$  such that
\begin{equation}
J(\Pi^*) =\sup_{\Pi\in\mathcal{A}}J(\Pi),
\end{equation}
where
\begin{equation}\label{eq7.5}
J(\Pi) := E[U(X^{\Pi}(T))].
\end{equation}
Note that, in terms of our process $X(t,z)$ we have
\begin{equation}\label{eq7.6}
\begin{cases}
dX(t,z) =b_0(t, t, z)\Pi(t)X(t,z)dt +\sigma_0(t, t, z)\pi(t,z)X(t,z)dB(t)\\
+(\int_0^t\frac{\partial b_0}{\partial t}(t,s,z)\pi(s,z)X(s,z)ds+\int_0^t\frac{\partial \sigma_0}{\partial t}(t,s,z)\pi(s,z)X(s,z)dB(s))dt; t\in[0, T] \\
X(0,z)=x_0>0
\end{cases}
\end{equation}
or in the integral form:
\begin{equation}\label{eq7.6b}
X(t,z) = x_0 + \int_0^t b_0(t, s, z)\pi(s,z)X(s,z)ds + \int_0^t \sigma_0(t, s,z)\pi(s,z)X(s,z)dB(s); t\in[0, T]
\end{equation}

We assume that $b_0(t, s,z) = b_0(t, s,z,\omega)$ and $\sigma_0(t, s,z) = \sigma_0(t, s,z, \omega)$ are given bounded processes, and that
$b_0(t, s)$ and $\sigma_0(t, s)$ are $\mathcal{F}_s$-measurable for all $s, t$ and $C^1$ with respect to
$t$ for all $s$, a.s.
We also assume that
\begin{equation}\label{eq7.6a}
\sigma_0(t, s,z) \geq c_0 \text{ a.s. for all t, } s\in[0, T] \text{ for some constant  }c_0 > 0.
\end{equation}
We assume that $x_0 > 0$. If $\pi\in \mathcal{A}$, then it follows that $X(t,z) > 0$ for all $t\in[0, T]$. To see this, note that
from \eqref{eq7.6} we get
\begin{align}\label{eq7.7}
X(t,z)= x_0\exp(\int_0^t\sigma_0(s, s,z)\pi(s,z)dB(s)
+ \int_0^t b_0(s, s,z)\pi(s,z) - \frac{1}{2}\sigma_0^2(s, s,z)\pi^2(s,z) +\alpha(s,z) ds)> 0,
\end{align}
where
\begin{align}
\alpha(s,z) :=\int_0^s\frac{\partial b_0}{\partial s}(s, r,z)\pi(r,z)X(r,z)dr +\int_0^s\frac{\partial\sigma_0}{\partial s}(s, r,z)\pi(r,z)X(r,z)dB(r).
\end{align}
The performance functional gets the form
\begin{equation}
J(\pi) = E[U(X(T,z))E[\delta_Z(z)|\mathcal{F}_T]]
\end{equation}
This is a problem of the type investigated in the previous sections (in the special case with
no jumps) and we can apply the results there to solve it, as follows:
The Hamiltonian gets the form, with $u = \pi$,
\begin{align}
&H(t, x, z, \pi, p, q)= b_0(t, t,z)\pi xp + \sigma_0(t, t,z)\pi xq\nonumber\\
&+ \int_t^T \frac{\partial b_0}{\partial s}(s, t,z)\pi xp(s,z)ds + \int_t^T \frac{\partial \sigma_0}{\partial s}(s, t,z)\pi xE[D_tp(s,z) |\mathcal{F}_t]ds
\end{align}
while the BSDE of the adjoint processes becomes
\begin{equation}\label{eq7.11}
\begin{cases}
dp  (t, z) =[b_0(t, t,z)\pi(t,z)p(t,z) + \sigma_0(t, t)\pi(t,z)q(t,z)\\
+ \int_t^T \frac{\partial b_0}{\partial s}(s, t,z)\pi(s,z)p(s,z)ds + \int_t^T \frac{\partial \sigma_0}{\partial s}(s, t,z)\pi(s,z)E[D_tp(s) |\mathcal{F}_t]ds]dt \\
+ q(t, z)dB(t); t\in[0, T]\\
p(T, z)=U'(X(T,z))E[\delta_Z(z)|\mathcal{F}_T]
\end{cases}
\end{equation}
Suppose there exists an optimal control $\pi\in\mathcal{A}$ for (5.18) with corresponding $\hat{X}, \hat{p}, \hat{q}$. Then,
\begin{equation}
E[\frac{\partial H}{\partial \pi}(t, \hat{X}(t,z), \pi, \hat{p}, \hat{q})|\mathcal{F}_t]_{\pi=\hat{\pi}(t)} = 0
\end{equation}

i.e.,
\begin{align}
&E[b_0(t, t,z)\hat{X}(t,z)\hat{p}(t,z) + \sigma_0(t, t,z)\hat{X}(t)\hat{q}(t,z)\nonumber\\
&+ \int_t^T \frac{\partial b_0}{\partial s}(s, t,z)\hat{X}(t,z)\hat{p}(s,z)ds + \int_t^T\frac{\partial \sigma_0}{\partial s}(s, t,z)\hat{X}(t,z)E[D_t\hat{p}(s,z)| \mathcal{F}_t]ds|\mathcal{F}_t] = 0.
\end{align}
Since $\hat{X}(t) > 0$, this is equivalent to
\begin{align}\label{eq7.14}
&b_0(t, t,z)\hat{p}(t) + \sigma_0(t, t,z)\hat{q}(t)\nonumber\\
&+E[\int_t^T \frac{\partial b_0}{\partial s}(s, t,z)\hat{p}(s)ds+ \int_t^T\frac{\partial \sigma_0}{\partial s}(s, t,z)E[D_t\hat{p}(s)|\mathcal{F}_t]ds|\mathcal{F}_t] = 0.
\end{align}
We deduce that the corresponding BSDE \eqref{eq7.11} reduces to
\begin{equation}
\begin{cases}
d\hat{p}(t,z) =\hat{q}(t,z)dB(t); 0\leq t \leq T\\
\hat{p}(T,z ) = U'(\hat{X}(T,z ))E[\delta_Z(z)|\mathcal{F}_T],
\end{cases}
\end{equation}
which has the unique solution
\begin{equation}\label{eq7.19}
\hat{p}(t,z) = E[U'(\hat{X}(T,z))E[\delta_Z(z)|\mathcal{F}_T]|\mathcal{F}_t], \quad \hat{q}(t,z) = D_t\hat{p}(t,z).
\end{equation}
Substituted \eqref{eq7.19} into \eqref{eq7.14}, this gives the equation
\begin{align}\label{eq7.17}
&E[b_0(t, t,z)U'(\hat{X}(T, z))E[\delta_Z(z)|\mathcal{F}_T] + \sigma_0(t, t,z)D_t\{U'(\hat{X}(T ,z))E[\delta_Z(z)|\mathcal{F}_T]\}\nonumber\\
&+ \int_t^T\frac{\partial b_0}{\partial s}(s, t,z)E[U'(\hat{X}(T ,z))E[\delta_Z(z)|\mathcal{F}_T]|\mathcal{F}_s]ds\nonumber\\
&+ \int_t^T\frac{\partial \sigma_0}{\partial s}(s, t,z)E[D_t\{U'(\hat{X}(T ,z))E[\delta_Z(z)|\mathcal{F}_T]\}|\mathcal{F}_s] ds |\mathcal{F}_t] = 0,
\end{align}
where we have used that
\begin{align}
D_tE[U'(\hat{X}(T ,z))E[\delta_Z(z)|\mathcal{F}_T]|\mathcal{F}_t] = E[D_t\{U'(\hat{X} (T ,z))E[\delta_Z(z)|\mathcal{F}_T]\}| \mathcal{F}_t] ,
\end{align}
Equation \eqref{eq7.17} can be simplified to
\begin{align}
&b_0(t, t,z)E[U'(\hat{X}(T,z ))E[\delta_Z(z)|\mathcal{F}_T]|\mathcal{F}_t]
+ \sigma_0(t, t,z)E[D_t\{U'(\hat{X}(T,z ))E[\delta_Z(z)|\mathcal{F}_T]\}|\mathcal{F}_t]\nonumber\\
&+E[\int_t^T\frac{\partial b_0}{\partial s}(s, t,z)U'(\hat{X}(T,z ))E[\delta_Z(z)|\mathcal{F}_T]ds|\mathcal{F}_t]\nonumber\\
&+E[\int_t^T\frac{\partial \sigma_0}{\partial s}(s, t,z)D_t\{U'(\hat{X}(T ,z))E[\delta_Z(z)|\mathcal{F}_T]\}ds|\mathcal{F}t] = 0,
\end{align}
or
\begin{equation}\label{eq7.21}
\sigma_0(T, t,z)D_tE[U'(\hat{X}(T ,z))E[\delta_Z(z)|\mathcal{F}_T]|\mathcal{F}_t] + b_0(T, t,z)E[U'(\hat{X}(T,z ))E[\delta_Z(z)|\mathcal{F}_T]|\mathcal{F}_t] = 0.
\end{equation}

By \eqref{eq7.6a} we see that \eqref{eq7.21} can be written
\begin{equation}
\frac{D_tY (t)}{Y (t)} = -\frac{b_0(T, t,z)}{\sigma_0(T, t,z)},
\end{equation}
where
\begin{equation}\label{eq7.23}
Y (t,z) = E[U'(\hat{X}(T,z))E[\delta_Z(z)|\mathcal{F}_T]|\mathcal{F}_t] .
\end{equation}
By the chain rule for Malliavin derivatives, we deduce from  that
\begin{equation}\label{eq7.24}
D_t(ln Y (t,z)) = - \frac{b_0(T, t,z)}{\sigma_0(T, t,z)}.
\end{equation}
On the other hand, since $Y (t,z)$ is a positive martingale, there exists an adapted process $\theta_0(t,z)$ such that
\begin{equation}
dY (t,z) = \theta_0(t,z)Y (t,z)dB(t)
\end{equation}
i.e.,
\begin{equation}\label{eq7.26}
Y (t,z) = Y (0,z)\exp(\int_0^t \theta_0(s,z)dB(s) - \frac{1}{2} \int_0^t \theta_0^2(s,z)ds).:
\end{equation}
From \eqref{eq7.26} we get
\begin{equation}\label{eq7.27}
D_t(\ln Y (t,z)) = D_t(\int_0^t \theta_0(s,z)dB(s) - \frac{1}{2}\int_0^t \theta_0^2(s,z)ds) = \theta_0(t,z),
\end{equation}
Comparing \eqref{eq7.24} and \eqref{eq7.27} we conclude that
\begin{equation}
\theta_0(t,z) = -\frac{b_0(T, t,z)}{\sigma_0(T, t,z)}
\end{equation}
and hence, by \eqref{eq7.23},
\begin{align}\label{eq7.29}
&E[U'(\hat{X}(T,z))E[\delta_Z(z)|\mathcal{F}_T]|\mathcal{F}_t] = Y (t,z)\nonumber\\
&= E[U'(\hat{X}(T,z))E[\delta_Z(z)|\mathcal{F}_T]]\exp(\int_0^t \theta_0(s,z)dB(s) - \frac{1}{2}\int_0^t \theta_0^2(s,z)ds).
\end{align}
It remains to find the constant
\begin{equation}
c= E[U'(\hat{X}(T,z))E[\delta_Z(z)|\mathcal{F}_T]] .
\end{equation}
From \eqref{eq7.29} with $t = T$ we get
\begin{equation}\label{eq7.31}
U'(\hat{X}(T,z))E[\delta_Z(z)|\mathcal{F}_T]= c\exp(\int_0^T\theta_0(s,z)dB(s) - \frac{1}{2}\int_0^T \theta_0^2(s,z)ds).
\end{equation}
To make this more explicit, we proceed as follows:
Define
\begin{equation}
M(t, z) := E[\delta_Z(z)|\mathcal{F}_t]
\end{equation}
Then by the generalized Clark-Ocone theorem
\begin{equation}
\begin{cases}
dM(t, z) = E[D_t\delta_Z(z)|\mathcal{F}_t]dB(t) = \Phi(t, z)M(t, z)dB(t)\\
M(0,z)=1
\end{cases}
\end{equation}
where
\begin{align}\label{eq7.36}
\Phi(t, z) = \frac{E[D_t\delta_Z(z)|\mathcal{F}_t]}{E[\delta_Z(z)|\mathcal{F}_t]}.
\end{align}
Solving this SDE for $M(t,z)$ we get
\begin{equation}
M(t,z) =\exp(\int_0^t \Phi(s, z)dB(s) - \frac{1}{2} \int_0^t \Phi^2(s, z)ds).
\end{equation}
Substituting this into \eqref{eq7.31} we get
\begin{equation}\label{eq7.38}
U'(\hat{X}(T,z))= c\exp(\int_0^T(\theta_0(s,z)-\Phi(s, z))dB(s) - \frac{1}{2}\int_0^T (\theta_0^2(s,z)-\Phi^2(s, z))ds).
\end{equation}
Then we can deduce that
\begin{equation}\label{eq7.39}
\hat{X}(T,z)=(U')^{-1} (c\exp(\int_0^T(\theta_0(s,z)-\Phi(s, z))dB(s) - \frac{1}{2}\int_0^T (\theta_0^2(s,z)-\Phi^2(s, z))ds))=F(c).
\end{equation}

On the other hand, if we define
\begin{equation}\label{eq7.40}
\hat{K}(t, s,z) := \sigma_0(t, s,z)\hat{\pi}(s,z)\hat{X}(s,z),
\end{equation}
then by \eqref{eq7.6b}, the pair $(\hat{X}, \hat{K})$ solves the following (Yong type) backward stochastic Volterra integral equation
(BSVIE)
\begin{equation}\label{eq7.41}
\hat{X}(t,z) = F(c) - \int_t^T \frac{b_0(t, s,z)} {\sigma_0(t, s,z)} \hat{K}(t, s,z)ds - \int_t^T \hat{K}(t, s,z)dB(s); 0\leq s\leq T.
\end{equation}
By Theorem 3.2 in [20] the solution of this equation is unique. Putting $t = 0$ and taking expectation in
\eqref{eq7.41}, we get
\begin{equation}\label{eq7.42}
x_0 = E[F(c)] - \int_0^T E[\frac{b_0(t, s,z)} {\sigma_0(t, s,z)} \hat{K}(t, s,z)]ds.
\end{equation}
This equation determines implicitly the value of $c$.
Hence by \eqref{eq7.39} we have found the optimal terminal wealth $\hat{X}(T,z)$. Then, finally we obtain the optimal
portfolio $\hat{\pi}$ by \eqref{eq7.40}.
Conversely, since the functions $x \rightarrow U(x)$ and $(x, \pi)\rightarrow H(t, x, \pi,\hat{p},\hat{q})$ are concave, we see that $\hat{\pi}$ found
above satisfies the conditions of Theorem 4.1, and hence $\hat{\pi}$ is indeed optimal.
We summarize what we have proved as follows:
\begin{theorem} Assume that $\sigma_0(t, s) > 0$ is bounded away from 0, for $s, t \in [0, T ]$.
Then, the optimal portfolio $\hat{\Pi}(t)$ for the problem \eqref{eq7.5} is
\begin{equation}
 \hat{\Pi}(t,Z)=\hat{\pi}(t,z)_{z=Z}=\int_{\mathbb{R}}\pi(z)\delta_Z(z)dz
 \end{equation}
 where
\begin{equation}
\hat{\pi}(s,z) =\frac{\hat{K}(t, s,z)}{\sigma_0(t, s,z)\hat{X}(s,z)}; s\in[0, T ],
\end{equation}
with $(\hat{X},\hat{K})$ is the unique solution of the BSVIE \eqref{eq7.41} with $F$ defined by \eqref{eq7.39}, and the constant $c$ is the
solution of \eqref{eq7.42}.
\end{theorem}

\textbf{Acknowledgment}\\
We are grateful to Bernt {\O}ksendal for helpful comments.

\end{document}